\documentclass[11pt]{amsart}
\usepackage{amsmath,amsthm,epsfig,amssymb}
\usepackage{xypic}
\input xy
\title{Lattice Platonic Solids and their Ehrhart polynomial}
\author{Eugen J. Ionascu}
\curraddr{Department of Mathematics\\ Columbus State University\\4225 University Avenue\\
Columbus, GA 31907\\
Honorific Member of the Romanian Institute of Mathematics ``Simion
Stoilow" } \email{ionascu@columbusstate.edu;} \subjclass{52C07, 05A15, 68R05}
\date{November $4^{th}$, 2011}
  \keywords{Ehrhart polynomial, linear Diophantine
equations, lattice points, cubes, tetrahedrons, octahedrons, Dedekind sums, sequences}
\begin{document}
\def\sms{\small\scshape}
\baselineskip18pt
\newtheorem{theorem}{\hspace{\parindent}
T{\scriptsize HEOREM}}[section]
\newtheorem{proposition}[theorem]
{\hspace{\parindent }P{\scriptsize ROPOSITION}}
\newtheorem{corollary}[theorem]
{\hspace{\parindent }C{\scriptsize OROLLARY}}
\newtheorem{lemma}[theorem]
{\hspace{\parindent }L{\scriptsize EMMA}}
\newtheorem{definition}[theorem]
{\hspace{\parindent }D{\scriptsize EFINITION}}
\newtheorem{problem}[theorem]
{\hspace{\parindent }P{\scriptsize ROBLEM}}
\newtheorem{conjecture}[theorem]
{\hspace{\parindent }C{\scriptsize ONJECTURE}}
\newtheorem{example}[theorem]
{\hspace{\parindent }E{\scriptsize XAMPLE}}
\newtheorem{remark}[theorem]
{\hspace{\parindent }R{\scriptsize EMARK}}
\renewcommand{\thetheorem}{\arabic{section}.\arabic{theorem}}
\renewcommand{\theenumi}{(\roman{enumi})}
\renewcommand{\labelenumi}{\theenumi}
\newcommand{\Q}{{\mathbb Q}}
\newcommand{\Z}{{\mathbb Z}}
\newcommand{\N}{{\mathbb N}}
\newcommand{\C}{{\mathbb C}}
\newcommand{\R}{{\mathbb R}}
\newcommand{\F}{{\mathbb F}}
\newcommand{\K}{{\mathbb K}}
\newcommand{\D}{{\mathbb D}}
\def\phi{\varphi}
\def\ra{\rightarrow}
\def\sd{\bigtriangledown}
\def\ac{\mathaccent94}
\def\wi{\sim}
\def\wt{\widetilde}
\def\bb#1{{\Bbb#1}}
\def\bs{\backslash}
\def\cal{\mathcal}
\def\ca#1{{\cal#1}}
\def\Bbb#1{\bf#1}
\def\blacksquare{{\ \vrule height7pt width7pt depth0pt}}
\def\bsq{\blacksquare}
\def\proof{\hspace{\parindent}{P{\scriptsize ROOF}}}
\def\pofthe{P{\scriptsize ROOF OF}
T{\scriptsize HEOREM}\  }
\def\pofle{\hspace{\parindent}P{\scriptsize ROOF OF}
L{\scriptsize EMMA}\  }
\def\pofcor{\hspace{\parindent}P{\scriptsize ROOF OF}
C{\scriptsize ROLLARY}\  }
\def\pofpro{\hspace{\parindent}P{\scriptsize ROOF OF}
P{\scriptsize ROPOSITION}\  }
\def\n{\noindent}
\def\wh{\widehat}
\def\eproof{$\hfill\bsq$\par}
\def\ds{\displaystyle}
\def\du{\overset{\text {\bf .}}{\cup}}
\def\Du{\overset{\text {\bf .}}{\bigcup}}
\def\b{$\blacklozenge$}

\def\eqtr{{\cal E}{\cal T}(\Z) }
\def\eproofi{\bsq}

\begin{abstract} First, we calculate the Ehrhart
polynomial associated to an arbitrary cube with integer coordinates for its vertices. Then, we use this result to derive relationships between the Ehrhart polynomials for regular lattice tetrahedrons and those for regular lattice octahedrons. These relations allow one to reduce the calculation of these polynomials to only one coefficient.
\end{abstract} \maketitle
\section{INTRODUCTION}
In the 1960's, $\rm Eug\grave{e}ne$ Ehrhart (\cite{ee},\cite{ee2}) proved that given a $d$-dimensional compact simplicial complex in $\mathbb
R^n$ ($1\le d\le n$), denoted here generically by $\cal P$,  whose vertices are in the lattice $\mathbb Z^n$,
there exists a polynomial  $L({\cal P}, t)\in \mathbb
\mathbb Q [t]$ of degree $d$, associated with $\cal P$, satisfying

\begin{equation} \label{ehrhartpolynomial}
L({\cal P},t)=the \ cardinality\  of\ \{t\cal P\}\cap \mathbb Z^n, \ t\in \mathbb N.
\end{equation}
It is known that

$$L({\cal P}, t) = Vol ({\cal P}) t^n + \frac{1}{2}Vol (\partial {\cal P})t ^{n-1} +... + \chi ({\cal P})$$

\n where  $\chi({\cal P})$ is the Euler characteristic of $\cal P$, and $Vol (\partial {\cal P})$  is the surface area of $\cal P$ normalized
with respect to the sublattice on each face of $\cal P$.

In \cite{ee0}, ${\rm Eug\grave{e}ne}$ Ehrhart has classified the regular convex polyhedra in $\mathbb Z^3$. It turns out that only cubes, regular tetrahedrons and regular octahedrons can be embedded in the usual integer lattice.  We arrived at the same result in \cite{ejips} using a construction of these polyhedrons from equilateral triangles, and this is how we got interested in this line of questioning.  This led us to the following simple description of all cubes in
$\mathbb Z^3$. One takes an odd positive integer, say $d$, and a primitive solution of the Diophantine equation $a^2+b^2+c^2=3d^2$
($\gcd(a,b,c)=1$). There are equilateral triangles in any plane having equation $ax+by+cz=f$, which can be parameterized  in terms of two integers $m$ and $n$.
The side-lengths of such a triangle are equal to $d\sqrt{2(m^2+mn+n^2)}$. In order to rise in space from such a triangle to form a regular tetrahedron, one needs to have satisfied  the necessary condition

\begin{equation}\label{carrollewiscondition}
m^2+mn+n^2=k^2\  for\ some\ odd \ k\in \mathbb N.
\end{equation}
Parenthetically speaking we remind the reader about Carroll Lewis ' conjecture on having infinitely many triples of Pythagorean triangles of equal area. The answer to that problem was to take solutions of (\ref{carrollewiscondition}) and construct triangles with sides $2uv$, $|u^2-v^2|$ and $u^2+v^2$ having $(u,v)\in\{ (m,k),(n,k),(m+n,k)\}$. This doesn't seem to be that strange of a coincidence, given the number of Diophantine equations that one need to solved in both cases. However, we could not find a reason of having (\ref{carrollewiscondition}) as a sufficient condition.

 If (\ref{carrollewiscondition}) is satisfied, there are two possibilities. If $k$ is a multiple of $3$, then one can complete the triangle
in both sides of the plane to a regular tetrahedron in
$\mathbb Z^3$, and if $k$ is not divisible by $3$, then one can complete the triangle in exactly one side to form a regular tetrahedron in $\mathbb Z^3$ (see Figure~1). Every such regular tetrahedron can then be completed to a cube in $\mathbb Z^3$ with side-lengths equal to $dk$. Every regular octahedron in $\mathbb Z^3$ is the dual of the double of a cube in $\mathbb Z^3$. We are going to make these constructions a little more specific in the last section.

\begin{figure}
\begin{center}\label{fig1}
\epsfig{file=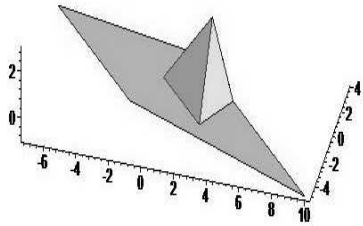,height=1.5in,width=2in}
 \end{center}\caption{Expand an equilateral triangle}
\end{figure}

It is natural to ask the question that we think  Ehrhart  himself asked: ``What is the form that
the polynomial in (\ref{ehrhartpolynomial}) takes for these regular lattice polyhedra?". The purpose of this paper is to answer this question for cubes (in a very simple way), and for regular tetrahedrons and octahedrons, the only such regular polyhedra in $\mathbb Z^3$.

For completeness and due credit, we include Ehrhart's  idea in \cite{ee3} to characterize all cubes in $\mathbb Z^3$.  This was based on a theorem of  Olinde  Rodrigues:  {\it The set of 3-by-3 orthogonal matrices can be given by four
real parameters $a$, $b$, $c$, $d$, not simultaneously zero, as follows:

\begin{equation}\label{olindeRodrigues}\frac{\pm 1}{a^2+b^2+c^2+d^2}\left[
                         \begin{array}{ccc}
                           a^2 + b^2- c^2-d^2 & 2(bc + da) & 2(bd-ca) \\
                           2(bc-da) & a^2- b^2+ c^2-d^2 & 2(cd+ ba) \\
                           2(bd + ca) & 2(cd-ba) & a^2 -b^2-c^2 +d^2\\
                         \end{array}
                       \right].
                       \end{equation}}

It is clear that every cube in $\mathbb Z^3$ can be translated so that a vertex becomes the origin and the three  vectors defined by the three sides starting from the origin give an orthogonal basis for $\mathbb R^3$. Hence, one can construct a 3-by-3 orthogonal matrix from these vectors which has rational entries. Conversely, one can construct a cube  in $\mathbb Z^3$ from such an orthogonal matrix which has rational entries. In what follows we will do this association so that the vectors (points) are determined by the rows.  The construction here is to take four integers $a$, $b$, $c$ and $d$ in (\ref{olindeRodrigues}), simplify by whatever is possible and then get rid of the denominators to obtain the three vectors with integer coordinates that determine the cube. This construction is similar to the classical parametrization of the Heronian triangles.

Our approach to the classification allows us to start in terms of the side lengths. However, Ehrhart's construction is useful in answering other questions about these objects. For instance, one can see that there are such cubes of any side length (other than the trivial ones, multiples of the unit cube) since every natural number can be written as a sum of four perfect squares. It turns out that there are only odd number side lengths for irreducible cubes, i.e. a cube which is not an integer multiple of a smaller cube in $\mathbb Z^3$.

Let us begin with some of the smallest irreducible cubes. We introduced then here by
orthogonal matrices with rational entries and defined up to the usual symmetries of the space (equivalent classes modulo the 48-order subgroup of all orthogonal matrices with entries $0$ or $\pm 1$, denoted by ${\cal S}_{o}$) . As we mentioned before, and this will make a difference, the cubes
are essentially determined by the rows. Obviously, the Ehrhart polynomials are identical  for all cubes in the same equivalence class.

We will denote the Ehrhart polynomial for an irreducible cube $C_{\ell}$ of side-length $\ell=2k-1$, $k\in \mathbb N$,   by $L(C_{\ell},t)$. From the general theory we must have

\begin{equation}\label{firsform}L(C_{\ell},t)=\ell^3 t^3+\lambda_1t^2+\lambda_2t+1,\ t\in \mathbb Z,
\end{equation}

\n where $\lambda_1$ is half the sum of the areas of the faces of the cube $C_{\ell}$, each face being normalized by the area of a fundamental domain of the sublattice contained in that face. The coefficient $\lambda_2$ is in general a problem (see, for example \cite{diazrobins1997}) but
in this case it takes a simple form as we will show in Section 3.

\n For the unit cube $C_1=I$ (the identity matrix), obviously, $L(C_1,t)=(t+1)^3$. There is only one cube (right or left equivalence classes modulo ${\cal S}_{o}$) for each $\ell=2k-1$ for $k=1,2,3,4,5$ and $6$: $C_1=I$,

$$C_3:= \frac{1}{3}
\left[ \begin{array}{rrr}
          -1 & 2 & 2 \\
          2 & -1 & 2 \\
          2 & 2 & -1 \\
        \end{array}
      \right], \ \ C_5:= \frac{1}{5}
\left[ \begin{array}{rrr}
          4 & 3 & 0 \\
          3 & -4 & 0 \\
          0 & 0& 5 \\
        \end{array}
      \right],\ \
      C_7:= \frac{1}{7}
\left[ \begin{array}{rrr}
          -2 & 6 & 3 \\
          3 & -2 & 6 \\
          6 & 3 & -2 \\
        \end{array}
      \right],
$$

$$
C_9:= \frac{1}{9}
\left[ \begin{array}{rrr}
          7 & 4 & -4 \\
          4 & 1 & 8 \\
          -4 & 8 & 1 \\
        \end{array}
      \right],\ \ and\ C_{11}:= \frac{1}{11}
\left[ \begin{array}{rrr}
          2 & 9 & 6 \\
          9 & 2 & -6 \\
          6 & -6 & 7 \\
        \end{array}
      \right].
      $$

\n For $k=7$ we have $C_{13}:= \frac{1}{13}
\left[ \begin{array}{rrr}
          -3 & 12 & 4 \\
          4 & -3 & 12 \\
          12 & 4 & -3 \\
        \end{array}
      \right]$, and an extra orthogonal matrix: $$\hat{C}_{13}:= \frac{1}{13}
\left[ \begin{array}{rrr}
          5 & 12 & 0 \\
          12 & -5 & 0 \\
          0 & 0 & 13\\
        \end{array}
      \right].$$ One peculiar thing about the Ehrhart polynomials associated with these cubes so far, is that there is an unexpected factor in their factorization:

$$\begin{array}{c}
L(C_3,t)=(3t+1)(9t^2+1),\   L(C_5,t)=(5t+1)(25t^2+2t+1), \ L(C_7,t)=(7t+1)(49t^2-4t+1), \\ \\  L(C_9,t)=(9t+1)(81t^2-6t+1),\ L(C_{11},t)=(11t+1)(121t^2-8t+1), \\ \\  L(C_{13},t)=(13t+1)(169t^2-10t+1),   \ and \ L({\hat C}_{13},t)=(13t+1)(169t^2+2t+1).
\end{array}
$$

\n This suggests that
\begin{equation}\label{conjecture}
L(C_{\ell},t)=(\ell t+1)(\ell^2t^2+\alpha t+1), \ t\in \mathbb Z, \ and\ some\
\alpha\in \mathbb Z.
\end{equation}

We can easily prove that this is indeed the case for cubes of a special form, like $C_5$ and $\hat{C}_{13}$ above. Let us consider a primitive Pythagorean triple
$(a,b,c)$, with $a^2+b^2=c^2$. In the $xy$-plane, we construct the square with vertices $O(0,0,0)$, $A(a,b,0)$, $B(a-b,a+b,0)$, and $C(-b,a,0)$.  We then translate this face along
the vector $c\overrightarrow{k}$ to form a cube of side-lengths equal to $c$. Let us denote this cube by $C_{a,b,c}$. It is easy to argue that because we have a primitive Pythagorean triple $(a,b,c)$, we have no lattice points on the sides of $OABC$, other than its vertices. The coefficient $\lambda_1$ in (\ref{firsform}), is equal to $\frac{1}{2}(c^2+c^2+4c)$ because two of the faces have to be normalized by $1$ and four of the faces have to be normalized by $c(\frac{c}{c})=c$. By Pick's theorem, applied to $OABC$,  we must have
$$c^2=\frac{\# \{points\ on\ the\ sides\}}{2}+\# \{interior\ points\ of\ OABC\}-1=\# \{interior\ points\ of\ OABC\}+1.$$
Hence the number of lattice points in the interior of $OABC$ is $c^2-1$. Therefore the number of lattice points in $C_{a,b,c}$ is $(c+1)(c^2+3)=c^3+c^2+3c+3$.
The polynomial $L(C_{a,b,c},t)=c^3t^3+(c^2+2c)t^2+(c+2)t+1=(ct+1)(c^2t^2+2t+1)$ satisfies exactly the condition $L(C_{a,b,c},1)=(c+1)(c^2+3)$. We then have shown that
(\ref{conjecture}) is true for infinitely many cubes $C_{\ell}$.

\begin{proposition} Given a primitive Pythagorean triple, $a^2+b^2=c^2$, the cubes in the class of   $C_{a,b,c}:=\frac{1}{c}
\left[ \begin{array}{rrr}
          a & b & 0 \\
          -b & a & 0 \\
          0 & 0 & c\\
        \end{array}
      \right]$, have the same  Ehrhart polynomial given by
$$L(C_c,t)=(ct+1)(c^2t^2+2t+1),\ t\in \mathbb N.$$
\end{proposition}

The general formula is proved in Section 3. Section 2 is basically dealing with the second coefficient in (\ref{firsform}). In Section 4, we look at the Ehrhart polynomial for regular tetrahedrons and  regular octahedrons with lattice vertices.  We show some nice relationship between the two and give a method of computing them in terms of Dedekind sums.

\section{The coefficient $\lambda_1$}

Let us prove the following lemma which is more or less contained in \cite{ejieeqtr} or it may already be a known result but we include it here for completeness.

\begin{lemma}\label{planefunddomain} Let $a$, $b$ and $c$ integers such that $gcd(a,b,c)=1$. The sub-lattice ${\cal L}:=\{(x,y,z)\in \mathbb Z^3| ax+by+zc=0\}$ is
generated by two vectors $\overrightarrow{u}$ and $\overrightarrow{\tau}$ such that $|\overrightarrow{u}\times \overrightarrow{\tau}|=\sqrt{a^2+b^2+c^2}.$
\end{lemma}
\n {\bf Proof:} The statement is clear if two of the numbers $a$, $b$ or $c$ are zero. Otherwise the vectors
$\overrightarrow{u}=\frac{1}{\gcd(a,b)}(-b,a,0)$,
$\overrightarrow{v}=\frac{1}{\gcd(a,c)}(-c,0,a)$ and
$\overrightarrow{w}=\frac{1}{\gcd(b,c)}(0,-c,b)$ correspond to
points in ${\cal L}$. We first show that
${\cal L}$ is generated by $\overrightarrow{u}$,
$\overrightarrow{v}$ and $\overrightarrow{w}$ (a point $P$ in
${\cal P}_{a,b,c}$ is identified by its position vector
$\overrightarrow{OP}$ as usual).  Let $(x,y,z)$ be an arbitrary point in $\cal L$. If we define $\omega:=gcd(a,b)$, then
because $gcd(a,b,c)=1$, we need to have $z=\omega z'$ with
$z'\in \mathbb Z$. Also, the existence of integers $k$ and
$l$ such that $ka+\ell b=\omega$ is insured by the fact
$\omega=gcd(a,b)$. This means that we have
$$(x,y,z)-z'[gcd(a,c)k\overrightarrow{v}+gcd(b,c)\ell \overrightarrow{w}]=(\alpha',\beta',0)\in {\cal L}.$$
Since $a\alpha'+b\beta'=0$ we see that $(\alpha',\beta',0)=\lambda
\overrightarrow{u}$ for some $\lambda\in \mathbb  Z^3$. This proves that the lattice $\cal L$ is generated by  $\overrightarrow{u}$,
$\overrightarrow{v}$ and $\overrightarrow{w}$. In fact, we showed even more: the the lattice $\cal L$ is generated by $\overrightarrow{u}$ and
$\overrightarrow{\tau}:=gcd(a,c)k\overrightarrow{v}+gcd(b,c)\ell
\overrightarrow{w}$, where $ka+\ell b=\omega$.

We observe that $\overrightarrow{u}\times \overrightarrow{v}=\frac{a}{\omega \gcd(a,c)}(a\overrightarrow{i}+b\overrightarrow{j}+c\overrightarrow{k})$
and similarly $\overrightarrow{u}\times \overrightarrow{w}=\frac{b}{\omega \gcd(b,c)}(a\overrightarrow{i}+b\overrightarrow{j}+c\overrightarrow{k})$.
Hence the area of the parallelogram determined by $\overrightarrow{u}$ and $\overrightarrow{\tau}$ is equal to

$$|\overrightarrow{u}\times \overrightarrow{\tau}|=|\frac{1}{\omega}(ak+bl)(a\overrightarrow{i}+b\overrightarrow{j}+c\overrightarrow{k})|=\sqrt{a^2+b^2+c^2}.\ \ \bsq$$

\n Let us now assume that we have an arbitrary cube in $ \mathbb Z^3$,

\begin{equation}\label{orthogonalmatrix}
C_{\ell}=\frac{1}{\ell}\left[ \begin{array}{rrr}
          a_1 & b_1 & c_1 \\
          a_2 & b_2 & c_2 \\
          a_3 & b_2 & c_3\\
        \end{array}
      \right],
\end{equation}

\n with $a_i$, $b_i$ and $c_i$ integers such that $a_ia_j+b_ib_j+c_ic_j=\delta_{i,j}\ell^2$ for all $i$, $j$ in $\{1,2,3\}$.
We define $d_i:=\gcd(a_i,b_i,c_i)$. It is clear that the $d_i$ are divisors of $\ell$. Let us also introduce the numbers $d'_i=\frac{\ell}{d_i}$, $i=1,2,3$.
Then, we have the following expression for the first coefficient in (\ref{firsform}).
\begin{theorem}\label{coefficeinta1}The coefficient $\lambda_1$ in (\ref{firsform}) is given by
\begin{equation}\label{coefficeinta1formula} \lambda_1=\ell (d_1+d_2+d_3).
\end{equation}
\end{theorem}

{\n Proof:} We use Lemma~\ref{planefunddomain} for each of the faces of the cube. Opposite faces contribute the same way into the calculation. Say we take a face containing the points $(a_1,b_1,c_1)$ and $(a_2,b_2,c_2)$. The irreducible normal vector to this face is clearly $\frac{1}{d_3}(a_3,b_3,c_3)$. The area of a fundamental domain here is given by
$\sqrt{\frac{1}{d_3^2}(a_3^2+b_3^2+c_3^2)}=d'_3$. By the general theory $\lambda_1=\frac{1}{2}(2\frac{\ell^2}{d'_1}+2\frac{\ell^2}{d'_2}+
 2\frac{\ell^2}{d'_3})=\ell(d_1+d_2+d_3)$.\eproof

 It is natural to ask at this point whether or not it is possible to have all of the $d_i$'s larger than one. It turns out that this is possible and as before, in our line of similar investigations, the first $\ell$ with this property is $\ell=1105=5(13)(17)$:

 \begin{equation}\label{firstexample} C_{1105}=\frac{1}{1105} \left[ \begin{array}{rrr}
          -65 & 156 & 1092 \\
          420 & 1015 & -120 \\
          1020 & -408 &119\\
        \end{array}
      \right].
\end{equation}

\begin{corollary}\label{columnversusrows}  For a matrix $C_{\ell}$ as in (\ref{orthogonalmatrix}), such that $C_{\ell}$ and $C_{\ell}^{-1}$ are in the same equivalence class (left or right) modulo ${\cal S}_{o}$, we have

 \begin{equation}\label{corollaryst}
d_1+d_2+d_3=\gcd(a_1,a_2,a_3)+\gcd(b_1,b_2,b_3)+\gcd(c_1,c_2,c_3).
\end{equation}
\end{corollary}

\n \proof.  The Ehrhart polynomial must be the same for the corresponding cubes in the same equivalence class. \eproof

\n We believe that this corollary applies to all $\ell<1105$, and of course to a lot of other cases, but we do not have a proof of this. A counterexample to the hypothesis of this corollary
is given by the matrix given in (\ref{firstexample}). In this case, $d_1+d_2+d_3=35$ and $\gcd(a_1,a_2,a_3)+\gcd(b_1,b_2,b_3)+\gcd(c_1,c_2,c_3)=7$.

\section{The coefficient $\lambda_2$}

The main idea in calculating the coefficient $\lambda_2$ is to take advantage of the fact that every cube defined by (\ref{orthogonalmatrix}) can be used to form a wandering set $W$ for the space under integer translations along the the vectors $\overrightarrow{\alpha}=(a_1,b_1,c_1)$, $\overrightarrow{\beta}=(a_2,b_2,c_2)$, and $\overrightarrow{\gamma}=(a_3,b_3,c_3)$, i.e.,

$$\mathbb R^3=\underset{i,j,k\in \mathbb Z} {\overset{\circ}{\bigcup}} (W+i\overrightarrow{\alpha}+j\overrightarrow{\beta}+k\overrightarrow{\gamma}),$$

\n where $ {\overset{\circ}{\bigcup}}$ means a union of mutually disjoint sets.

\begin{figure}
\begin{center}
\epsfig{file=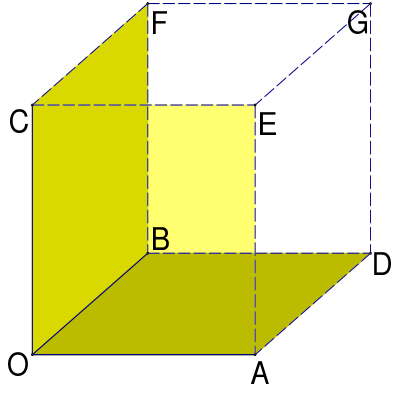,height=2in,width=2in}
 \end{center} \caption{Wandering set determined by the cube}\label{fig:fig22}
\end{figure}

The wandering set $W$ we will consider here, associated with a generic cube as in Figure~\ref{fig:fig22}, is the set of all points formed by the interior points of the cube to which we add the points of the faces OAEC, OADB and OBFC except the (closed) edges AD, DB, BF, FC, CE, and EA.  It is easy to see that such a set is indeed a wandering set. In our setting we think of $\overrightarrow{\alpha}$, $\overrightarrow{\beta}$ and $\overrightarrow{\gamma}$ as the vectors $\overrightarrow{OA}$, $\overrightarrow{OB}$ and $\overrightarrow{OC}$.
With the notation from the previous section, we have the following result.

\begin{theorem}\label{lastcoefficient}
The coefficient $\lambda_2$ in (\ref{firsform}) is equal to $d_1+d_2+d_3$.
\end{theorem}
\n \proof. Let us denote by $k$ the number of lattice points in $W$. For $n\in \mathbb N$, then the number of lattice points in
$$\underset{i,j,k\in \{1,2,...,n\}} {\overset{\circ}{\bigcup}} (W+i\overrightarrow{\alpha}+j\overrightarrow{\beta}+k\overrightarrow{\gamma}),$$

\n is equal to $n^3k$. On the other hand this number  is equal to $L(C_{\ell},n)+K$, where $K$ is the number of lattice points on three big faces of $nC$. It is easy to see that $K$ is $O(n^2)$ and so $k=\lim_{n\to \infty}  \frac{1}{n^3}(L(C_{\ell},n)+O(n^2))=\ell^3$.
We now need the following result which is well known.

\begin{theorem}\label{em}  ({\bf Ehrhart-Macdonald Reciprocity Law}) Given a compact simplicial lattice complex  $P$ (as before)  of dimension $n$, then
$$L(\overset{\circ}{ P},t) = (-1)^nL(P,-t), \ t\in \mathbb N.$$
\end{theorem}

\n Hence, according to this theorem the number of lattice points in the interior of $C_{\ell}$ is $\ell^3-\lambda_1+\lambda_2-1$. So, the number of lattice points on the boundary of $C_{\ell}$ is $2\lambda_1+2$. Let us denote by $\sigma$  the number of lattice points in the interiors of the sides  $\overline{OA}$, $\overline{OB}$ and $\overline{OC}$. Then we must have

$$2\lambda_1+2=2[k-(\ell^3-\lambda_1+\lambda_2-1)]+2\sigma+6 \ \Rightarrow \ \lambda_2 =\sigma+3.$$
Since, $\sigma =(d_1-1)+(d_2-1)+(d_3-1)$ the claim follows.\eproof

Putting these facts together we obtain.

\begin{theorem}\label{ehrhartpforcubes}
Given a cube $C_{\ell}$ constructed from a matrix as in (\ref{orthogonalmatrix}), its Ehrhart polynomial is given by
\begin{equation}\label{ep}
L(C_{\ell},t)=(\ell t+1)[\ell^2 t^2+(d_1+d_2+d_3-\ell)t+1], \ t\in \mathbb Z.
\end{equation}
\end{theorem}

There are some natural questions at this point. One of them is: ``What is the maximum number of lattice points that one can contain in a lattice cube of side lengths $\ell$?" We have the following corollary to the above theorem.

\begin{corollary}\label{corehrhartpforcubes}
Given a cube $C_{\ell}$ constructed from a matrix as in (\ref{orthogonalmatrix}) , the
maximum of lattice points inside or on the boundary of this cube cannot be more than  $(\ell+1)^3$. This value is attained for the cube $\ell C_1$.
\end{corollary}

\n \proof. Since $d_i$ is a divisor of $\ell$ we must have $d_i\le \ell$ and so the corollary follows from (\ref{ep}).\eproof

\n What  is the  answer to above question, if we restrict the problem to irreducible cubes? It is clearly a more complicated problem and it depends heavily on
$\ell$.

Finally, another question which seems natural here is: ``What is the cardinality of the set of all Ehrhart polynomials associated to irreducible cubes of sides $\ell$?"

\section{Regular Tetrahedrons and regular Octahedrons}

\begin{figure}
\begin{center}
\epsfig{file=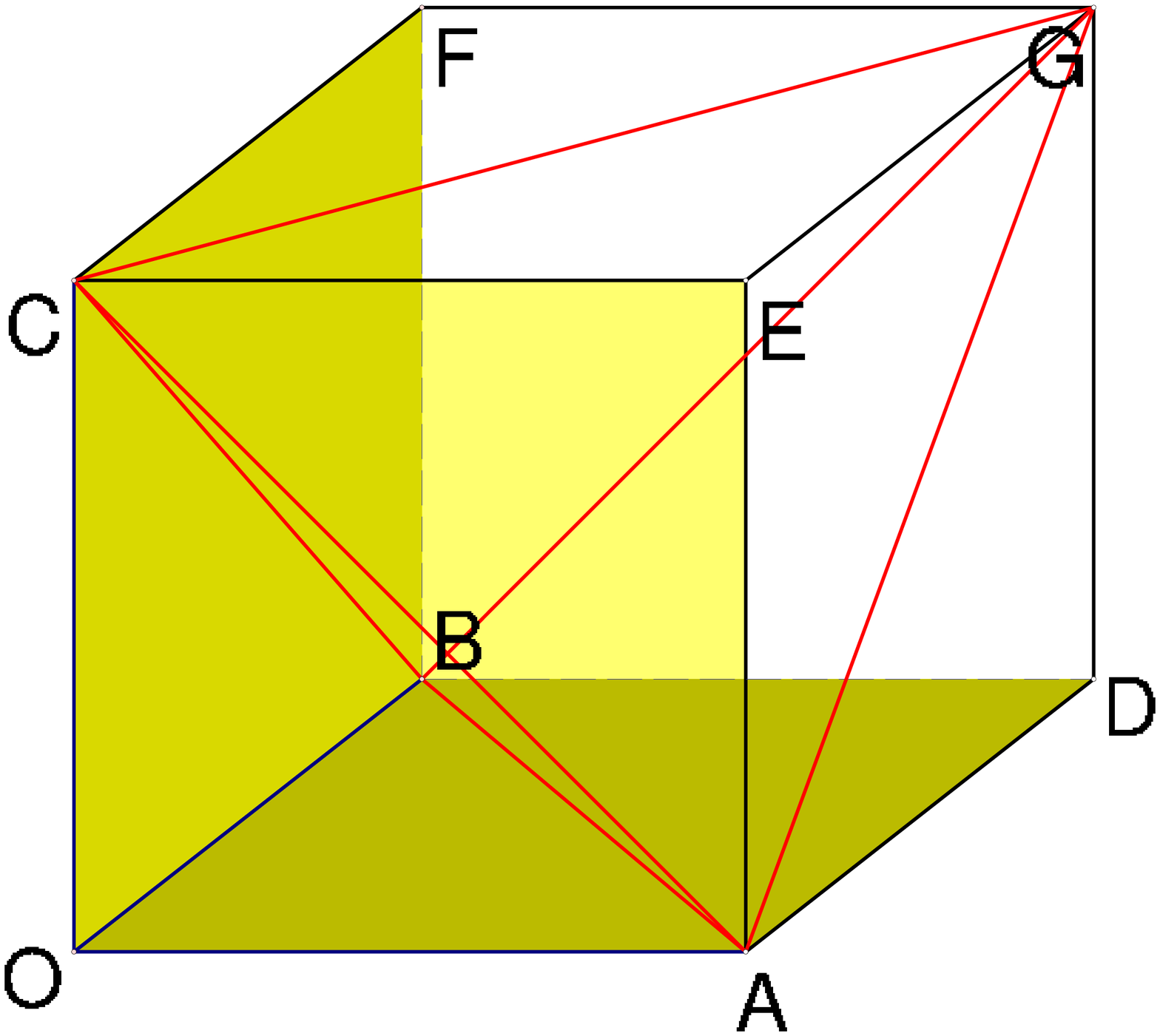,height=2in,width=2in}\epsfig{file=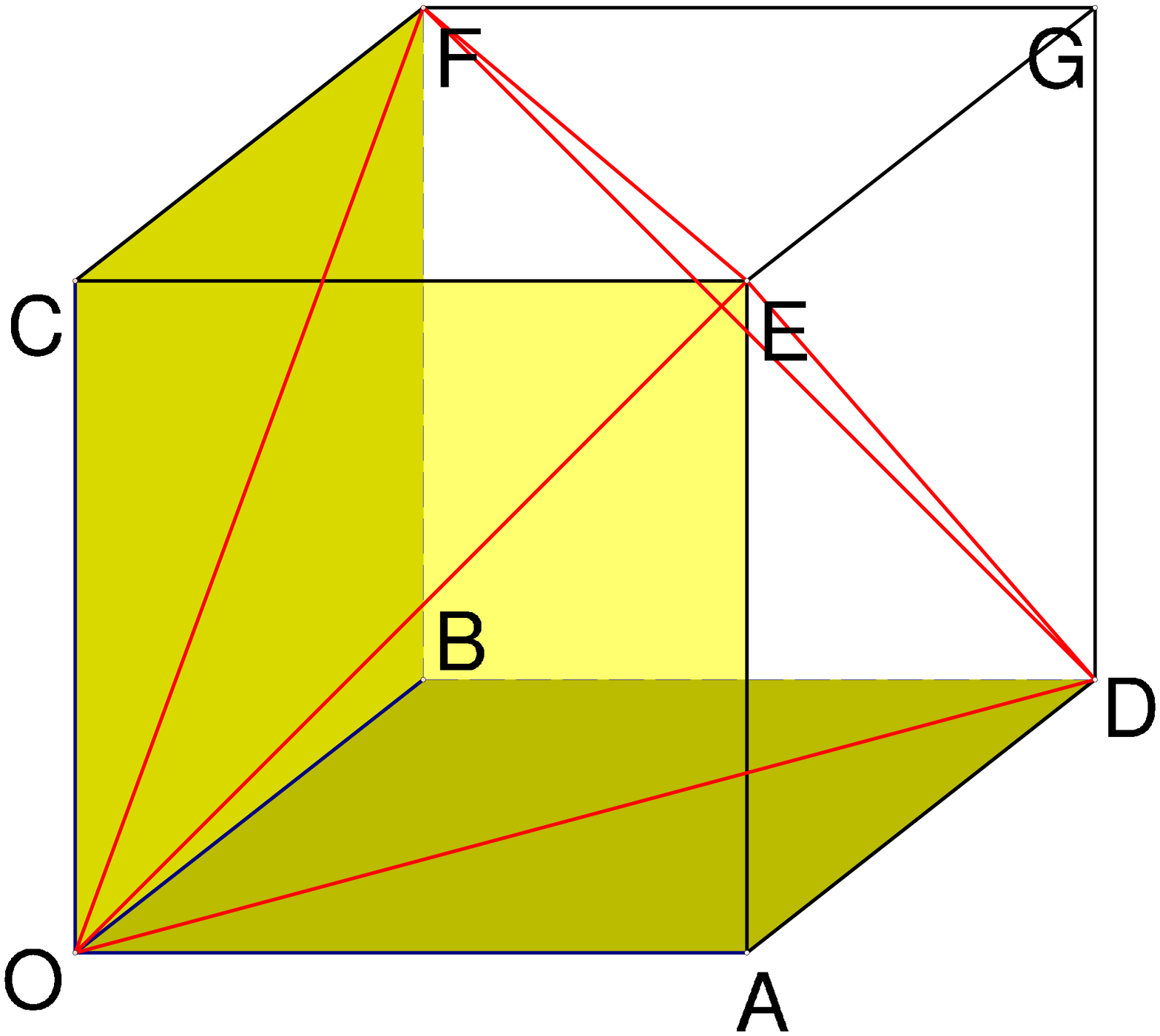,height=2in,width=2in}
 \end{center}\caption{A cube determines two tetrahedrons}\label{fig:tetrahedronsincubes}
\end{figure}

 We remind the reader that the cube in space (Figure~\ref{fig:tetrahedronsincubes}) is determined by an orthogonal  matrix as in (\ref{orthogonalmatrix}) by taking its vertices $O$ (the origin), $A$, $B$, $C$, $D$, $E$, $F$ and $G$, whose position vectors are  $\overrightarrow{OA}=\overrightarrow{\alpha}=(a_1,b_1,c_1)$, $\overrightarrow{OB}=\overrightarrow{\beta}=(a_2,b_2,c_2)$, $\overrightarrow{OC}=\overrightarrow{\gamma}=(a_3,b_3,c_3)$, $\overrightarrow{OD}=\overrightarrow{\alpha}+\overrightarrow{\beta}$, $\overrightarrow{OF}=\overrightarrow{\beta}+\overrightarrow{\gamma}$, $\overrightarrow{OE}=\overrightarrow{\gamma}+\overrightarrow{\alpha}$,  and $\overrightarrow{OG}=\overrightarrow{\alpha}+\overrightarrow{\beta}+\overrightarrow{\gamma}$.

In \cite{ejips}, we rediscovered Ehrhart's characterization (\cite{ee0}) of all regular polyhedra which can be imbedded in
$\mathbb Z^3$. Only cubes, tetrahedrons and octahedrons exist in
$\mathbb Z^3$ and there are infinitely many in each class. We have constructed all of these out of equilateral triangles.
In general, once a tetrahedron is constructed, it can always be completed to a cube. Vice versa, for a cube given by
(\ref{orthogonalmatrix}), there are two regular tetrahedrons which are shown in Figure~\ref{fig:tetrahedronsincubes}. They are in the same equivalence class modulo
the orthogonal matrices with entries $\pm 1$, denoted earlier by ${\cal S}_0$. Regular octahedrons can be obtained by doubling the coordinates of the cube $C_{\ell}$ and then taking the
centers of each face. This procedure is exhaustive. An octahedron in the same class can be obtained by simply taking the vertices whose position vectors  are $\pm \overrightarrow{\alpha}$,
$\pm \overrightarrow{\beta}$, and $\pm \overrightarrow{\gamma}$. We will going to use the notations $T_{\ell}$ and $O_{\ell}$ for the tetrahedrons and octahedrons
constructed in this way from $C_{\ell}$. If we deal with irreducible such objects, i.e. they cannot be divided by an  integer to obtain a strictly smaller similar polyhedron, then we may assume that $\ell$ is odd. The $T_{\ell}$ and $O_{\ell}$ have side-lengths equal to $\ell\sqrt{2}$. From the general Ehrhart theory (see \cite{beckAndRobins2007textbook}) we must have

\begin{equation}\label{firstfromto}
L(T_{\ell},t)=\frac{\ell^3}{3}t^3+\mu_1t^2+\mu_2t+1,\ \ L(O_{\ell},t)=\frac{4\ell^3}{3}t^3+\nu_1t^2+\nu_2t+1.
\end{equation}

\n Let us first look at some of the examples with the smallest side-lengths.
$$T_1:=\left[ \begin{array}{rrr}
          1 & 1 & 0 \\
          0 & 1 &  1 \\
          1 & 0 & 1 \\
        \end{array}
      \right],\ O_1:=
\left[ \begin{array}{rrr}
          1 & 0 & 0 \\
          0 & 1 & 0 \\
          0 & 0 & 1 \\
        \end{array}
      \right], \ with $$

      $$L(T_1,t)=\frac{t^3}{3}+t^2+\frac{5t}{3}+1,\ and \ L(O_1,t)=\frac{4}{3}t^3+2t^2+\frac{8t}{3}+1.$$
For the next sides,
$$T_3:=\left[ \begin{array}{rrr}
          1 & 1 & 4 \\
          1 & 4 &  1 \\
          4 & 1 & 1 \\
        \end{array}
      \right],\ O_3:=
\left[ \begin{array}{rrr}
          -1 & 2 & 2 \\
          2 & -1 & 2 \\
          2 & 2 & -1 \\
        \end{array}
      \right], \ with $$
   $$L(T_3,t)=9t^3+\frac{9}{2}t^2+\frac{13t}{2}+1,\ and \ L(O_3,t)=36t^3+9t^2-t+1.$$

\n These polynomials were computed initially with the help of a computer.

\subsection{Coefficients $\mu_1$ and $\nu_1$}

From the general theory we know that these coefficients can be computed in terms of the areas of their faces normalized by the area of the fundamental domains of the
sub-lattice of $\mathbb Z^3$ corresponding to that face.  Since the number of faces for $O_{\ell}$ is twice as big as for $T_{\ell}$ and since the parallel faces are
in the same equivalence class, we must have $\nu_1=2\mu_1$.

Let us introduce the divisors

\begin{equation}\label{divisors}
\begin{array}{l}
D_1=\gcd(a_1+a_2+a_3,b_1+b_2+b_3,c_1+c_2+c_3),\ \\ \\  D_2=\gcd(a_1+a_2-a_3,b_1+b_2-b_3,c_1+c_2-c_3), \\ \\ D_3=\gcd(a_1-a_2+a_3,b_1-b_2+b_3,c_1-c_2+c_3),\ and \ \ \\ \\  D_4=\gcd(-a_1+a_2+a_3,-b_1+b_2+b_3,-c_1+c_2+c_3).
\end{array}
\end{equation}
\n Observe that the vectors $\overrightarrow{\alpha}+\overrightarrow{\beta}+\overrightarrow{\gamma}$, $\overrightarrow{\alpha}+\overrightarrow{\beta}-\overrightarrow{\gamma}$, $\overrightarrow{\alpha}- \overrightarrow{\beta}+\overrightarrow{\gamma}$, $-\overrightarrow{\alpha}+\overrightarrow{\beta}+\overrightarrow{\gamma}$ are vectors normal to the faces of the $T_{\ell}$. By Lemma~\ref{planefunddomain},
we see that the area of each fundamental domain corresponding to a face of $T_{\ell}$ is given by one of the numbers $\frac{\ell \sqrt{3}}{D_i}$.

\begin{proposition} The coefficients $\mu_1$ and $\nu_1$ in (\ref{firstfromto}), are given by  \begin{equation}\label{a1b1}\mu_1=\frac{\nu_1}{2}=\ell(D_1+D_2+D_3+D_4)/4.\end{equation}
\end{proposition}

\n This explains the next examples which were obtained by brute force computations with Maple:

$$T_5:=
\left[ \begin{array}{rrr}
          7 & -1 & 0 \\
          4 & 3 & 5 \\
          3 & -4& 5 \\
        \end{array}
      \right],\ \
      O_5:=
\left[ \begin{array}{rrr}
          4 & 3 & 0 \\
          3 & -4 & 0 \\
          0 & 0 &5 \\
        \end{array}
      \right], \ with $$

 $$L(T_5,t)=\frac{125}{3}t^3+5t^2+\frac{1}{3}t+1,\ and \ L(O_5,t)=\frac{500}{3}t^3+10t^2+\frac{16}{3}t+1.$$

\subsection{Coefficients $\mu_2$ and $\nu_2$}

Let us observe that the cube in Figure~\ref{fig:tetrahedronsincubes} can be decomposed into four triangular pyramids OABC, DABG, FCGB, and EGCA, which can be
translated and some reflected through the origin to form half of $O_{\ell}$, and the regular tetrahedron $T_{\ell}$. We remind the reader of some notation we used in
the proof of Theorem~\ref{lastcoefficient}: we denoted by $\sigma$ the number of lattice points on the interiors of the edges $\overline{OA}$, $\overline{OB}$, and $\overline{OC}$. We showed that $\sigma=d_1+d_2+d_3-3$.

Let us balance the number $M$ of interior lattice  points
of $C_{\ell}$ using the above decomposition. According to Theorem~\ref{ehrhartpforcubes} and Theorem~\ref{em} we have $$M=-L(C_{\ell},-1)=(\ell-1)(\ell^2-d_1-d_2-d_3+\ell+1)=\ell^3-(d_1+d_2+d_3)(\ell-1)-1$$
 Some of the lattice points counted in $M$ are in the regular tetrahedron, and are counted by $L(T_{\ell},1)=\frac{\ell^3}{3}+\mu_1+\mu_2+1$. From  these we need to subtract the number of lattice points on the interiors of its sides, which we will denote by $\tau$, and also subtract 4 for its vertices. The rest of the points counted in $M$, are in the interiors of the four pyramids. If we multiply this number by two and add the number of lattice points in the interiors of the cube faces less $\tau$, we get the number of interior points of $O_{\ell}$ minus $2\sigma +1$. The number of interior lattice points of the cube faces is equal to $2\lambda_1+2-4\sigma-8$. In other words, we have

$$2(M-\frac{\ell^3}{3}-\mu_1-\mu_2-1+\tau+4)+2\lambda_1+2-4\sigma-8-\tau=\frac{4\ell^3}{3}-\nu_1+\nu_2-1-(2\sigma+1).$$

Taking into  account that $\nu_1=2\mu_1$ and $\lambda_1=(d_1+d_2+d_3)\ell=(\sigma+3)\ell$, this can be simplified to

\begin{equation}\label{cubeinoct}
\nu_2+2\mu_2=6+\tau.
\end{equation}

We close this section concluding what we have shown.

\begin{theorem}\label{sumofthelasttwocoefficeints} For a regular tetrahedron $T_{\ell}$ and  a regular octahedron $O_{\ell}$ constructed as before from an orthogonal matrix
with rational coefficients as in (\ref{orthogonalmatrix}), the coefficients $\mu_2$ and $\nu_2$ in (\ref{firstfromto}) satisfy

\begin{equation}\label{thesumproperty}
\nu_2+2\mu_2=(d_1+d_2+d_3+d_4+d_5+d_6),
\end{equation}

\n where $d_1$, $d_2$, $d_3$ are defined as before and $d_4=\gcd(a_1-a_2,b_1-b_2,c_1-c_2)$, $d_5=\gcd(a_1-a_3,b_1-b_3,c_1-c_3)$, and
$d_6=\gcd(a_3-a_2,b_3-b_2,c_3-c_2)$.
\end{theorem}

We have tried to find another relation that will help us find the two coefficients but it seems like there is not an easy way to
avoid, what are called in \cite{beckAndRobins2007textbook}, the building blocks of the lattice-point enumeration, the Dedekind sums.
This numbers require a little more computational power and we are wondering if a shortcut doesn't really exist. One would expect that the answer to our questions, for such regular objects, is encoded in the coordinates of their vertices in a relatively simple way.

\end{document}